\documentclass[preview]{elsarticle}

\usepackage{lineno,hyperref}
\usepackage{commath,amsmath, amssymb,float,caption}
\modulolinenumbers[5]

\journal{Journal of \LaTeX\ Templates}

 \newtheorem{thm}{Theorem}[section]
 \newproof{pf}{Proof}

\begin{document}

\begin{frontmatter}
\title{On fractional annealing process}

\author{Tran Manh Tuong}
\address{University of Finance and Marketing, Ho Chi Minh City, VietNam\\ email:tmtuong@ufm.edu.vn}

\author{Tran Hung Thao   \fnref{tmt}}
\address{Institute of Mathematics, Vietnam Academy of Science and Technology\\ thaoth@math.ac.vn}

\cortext[tmt]{Corresponding author}

\begin{abstract}

The paper deals with a full annealing process, for example for a metal alloy perturbed by some fractional noise caused by conditions in furnace such as initial temperature, initial structure of material itself, \dots The energy of the heating matter system depending on its states at each time, we consider some steady state at full annealing range by introducing a fractional model for it. And we determine states of the system related to this steady state.
\end{abstract}


\end{frontmatter}

\linenumbers

\section{Introduction}
As we know, full anneal  is a  special stage of annealing process where a condensed matter (a steel for example) is heated to slightly above the critical temperature (name austenitic temperature) and maintaining the temperature for a special period of time to allow the material  to fully form austenite or austenite cementite grain structure (see [1], [2], [4], [6]).

Denote by $X_t (t \ge 0), X_t \in \mathbb{R}^d$, a stochastic process presenting the state of the material under full annealing stage at time $t$ and $g(X_t)$ is the system energy corresponding to state $X_t$, where $g$ is a Borel function on $\mathbb{R}^n$.

It is well-know that  a classical stochastic model for annealing states $X_t$  is given by the following equation:
\begin{equation}\label{e1.1}dX_t = - \bigtriangledown g(X_t)dt + \sqrt{2T} dW_t  \tag{1.1}\end{equation}
where $W_t$ is a $d - $dimensional Brownian motion expressing a \emph{no-memory} noise and $T$ is the temperature of the system (see [3]). In our  context $T$ can be considered unchange as it is the austenitic temperature.

A problem is that full annealing treatments should be made in batch-type furnaces that are in full compliance of temperature uniformity and accuracy requirements of pyrometry AMS2750. In such type of furnace, the system state $(X_t, t \ge 0)$ is no more no-memory process because each value of $X_t$ at $t$ can influence upon its values some times after $t$. Therefore the model (1.1) is no more suitable. And we propose the following model
\begin{equation} \label{e1.2} dX_t = - \bigtriangledown g(X_t )dt + \sqrt{2T}dB^H_t \tag{1.2}\end{equation}
where: 
\begin{itemize}
	\item $B^H_t$ is a $d-$dimensional fractional Brownian motion of Liouvill form
		\[ \begin{split} B^H_t & = (B^{H_j}_t, j = 1,2,\dots,d), \\
				       B^{H_j}_t& = \int_0^t(t-s)^{H_j-1/2}dW^j_s, j = 1,2, \dots, d	 \end{split}\]
	\item $H_j$ are Hurst parameters, $0 < H_j < 1$.
	\item $W^j_t$ are standard Brownian motion.
	\item $g = g(x_1, \dots, x_d)$ is a Borel function $\mathbb{R}^d \to \mathbb{R}$ and $g \in C^1(\mathbb{R}^d)$.
\end{itemize}

The state process $(X_t)$ satisfying the equation (1.2) driven by the perturbation $B^H_t$ is a process with memory as expected.

But (1.2) is a fractional stochastic differential equation and its solution can not be found by means of traditional stochastic calculus as we can made for the  equation (1.2). Many approaches have been introduced to overcome this difficulty for such kind of equation. And one of these approaches is given by Thao T.H via an approximation method. This method is presented briefly in Section 2 below.

Now let $X^* = (X^*_1, X^*_2, \dots, X^*_d)$ be a steady state on the full annealing range that means $\frac{\partial g}{\partial x_i}(X^*) = 0, i = 1,2, \dots, d$. In this note we try  to find the relation between a state process satisfying (1.2) and the steady state $X^*$.

\section{Fractional Brownian motion and an approximation approach}
Now we recall some facts of one-dimensional fractional Brownian motion in Mandelbrot form ([7]-[10]). It is  a centered Gaussian process $W^H_t$ such that its covariance function is given by
\begin{equation} \label{e2.1} R(t, s) = \frac{1}{2}(t^{2H} + s^{2H}-\abs{t-s}^{2H}), 0 < H < 1 \tag{2.1}\end{equation}
If  $H = 1/2, R(t, s) = \min(t, s)$ and $W^H_t$ is an usual standard Brownian motion.\\
If $H \ne 1/2, W^H_t$ is a process of long memory, and Ito calculus cannot applied to it. And an approximation approach is given as follows ([7]).

It is known that $W^H_t$ can be decomposed as
\[ \label{e2.2} W^H_t = \frac{1}{\Gamma(\alpha+1)} \left[Z_t +  \int_0^t(t-s)^\alpha dW_s \right], \alpha = H - 1/2 \tag{2.2}\]
where $Z_t$  is a process having absolutely continuous paths and the long memory property focus at the term
\[ B^H_t = \int_0^t (t-s)^\alpha dW_s, \alpha = H - 1/2 \tag{2.3}\]
which is called the fractional Brownian motion of the Liouville form.
Since $B^H_t$ is not a semimartingale we introduce a new process for each $\epsilon > 0$
\[ \label{e2.4} B^{H,\epsilon}_t = \int_0^t (t-s + \epsilon)^\alpha dW_s \tag{2.4}\]
And the following two important facts have been proved in [7] and [10].
\begin{enumerate}[(i)]
	\item $B^{H,\epsilon}_t$ is a semimartingale:
		\[ dB^{H,\epsilon} = \alpha\varphi^\epsilon_tdt + \epsilon^\alpha dW_t \]
	where
		\[ \varphi_t^\epsilon = \int_0^t(t-s + \epsilon)^{\alpha-1}dW_s  \]
	\item $B^{H,\epsilon}_t$ converges to $B^H_t$ in $L^2(\Omega)$ as $\epsilon$ tends to 0. And this convergence is uniform with respect to $t$ belong to any finite interval $I$ of time, and we have also
	\[ \label{e2.5} \sup_{t\in I} \norm{B^{H,\epsilon}_t - B^H_t} \le C(\alpha)\epsilon^{2H} \tag{2.5}\]
	where $C(\alpha)$ is a positive constant depending only to $\alpha = H - 1/2$.
\end{enumerate}
Apart from this, a study of fractional stochastic differential equations is presented  in [8], [9] and [10].

\begin{figure}[H]
	\centering
  	\includegraphics[width=250pt]{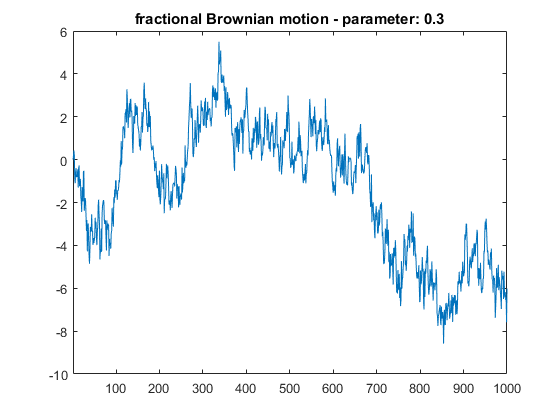}
  	\caption*{Fractional Brownian motion with H =0.3 and d = 1}
\end{figure}
\begin{figure}[H]
	\centering
	\includegraphics[width=250pt]{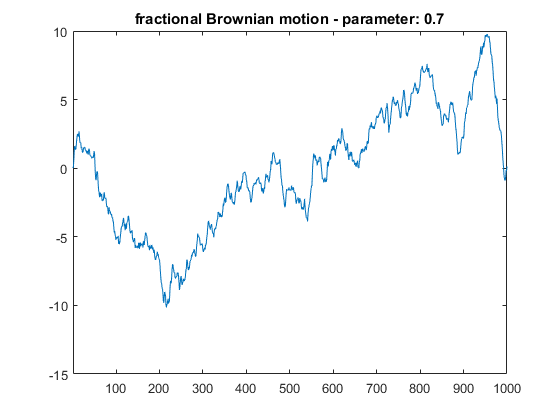}
  	\caption*{Fractional Brownian motion with H =0.7 and d = 1}
\end{figure}

\section{Fractional annealing process}
For the sake of the simplicity we consider the case of two dimensions $d = 2$. The case of general dimension $d$ can be naturally extended. So we are supposed to consider the following model
\[ \label{e3.0} dX_t = - \bigtriangledown g(X_t)dt + \sqrt{2T}dB^H_t \tag{3.0}\]
where
\[ X_t =  \begin{pmatrix} X^1_t \\ X^2_t \end{pmatrix}, g = g(x_1,x_2), B^H_t = \begin{pmatrix} B^{H_1}_t \\ B^{H_2}_t \end{pmatrix}\]
\[B^{H_j}_t = \int_0^t (t-s)^{\alpha_j}dW_j, \alpha_j = H_j -1/2, 0 < H_j < 1, j= 1,2 \]

Let $X^* = (X_1^*, X_2^*)$ be a steady state on full annealing range where $\frac{\partial g}{\partial x_1}(X^*) = 0, \frac{\partial g}{\partial x_2}(X^*) = 0$. We have to calculate state satisfying (1.6) and relating to $X^*$.

\subsection{A change of variables}
Denote by $F_1(x_1,x_2)$ and $F_2(x_1,x_2)$ be the partial derivatives of $g(x_1,x_2)$:
\[ \begin{matrix} F_1(x) = F_1(x_1,x_2) = -\frac{\partial}{\partial x_1}g(x_1, x_2) \\ F_2(x) = F_2(x_1,x_2) = -\frac{\partial}{\partial x_2}g(x_1, x_2) \end{matrix}  \]
and
\[ F(x) = \begin{pmatrix} F_1(x) \\ F_2(x) \end{pmatrix} \]

The equation (1.2) can be rewritten as following
\[ \label{e3.1.1} dX_t = F(X_t)dt + \sqrt{2T}dB^H_t \tag{3.1.1}\]
or
\[\begin{cases} dX_t^1 = F_1(X_t)dt + \sqrt{2T}dB^{H_1}_t\\ dX_t^2 = F_2(X_t)dt + \sqrt{2T}dB^{H_2}_t \end{cases}\]
We have now
\[ F_1(X^*) = 0 \text { and }F_2(X^*) = 0\]
Define a change of variable
\[ \label{e3.1.2} X_t = \begin{pmatrix} X^1_t \\ X^2_t \end{pmatrix} \to U_t = \begin{pmatrix} U^1_t \\ U^2_t \end{pmatrix} = X_t - X^* = \begin{pmatrix} X^1_t  - X^{1*}\\ X^2_t - X^{2*} \end{pmatrix} \tag{3.1.2}\]
We can see that
\[\label{3.1.3} \begin{pmatrix} dU^1_t \\ dU^2_t \end{pmatrix} = 
\begin{pmatrix} 
	\frac{\partial F_1}{\partial x_1}(X^*) & \frac{\partial F_1}{\partial x_2}(X^*)\\ 
	\frac{\partial F_1}{\partial x_2}(X^*) & \frac{\partial F_2}{\partial x_2}(X^*)
\end{pmatrix}
\begin{pmatrix} U^1_t \\ U^2_t \end{pmatrix}dt  + \sqrt{2T} 
\begin{pmatrix} dB^{H_1}_t \\ dB^{H_2}_t \end{pmatrix} \tag{3.1.3}\]
or
\[ \label{e3.1.4} dU_t = AU_tdt + \sqrt{2T}dB^H_t \tag{3.1.4}\]
where
\[ \label{e3.1.5} A= \begin{pmatrix} a_1&b_1 \\a_2&b_2 \end{pmatrix}: =  \begin{pmatrix} 
	\frac{\partial F_1}{\partial x_1}(X^*) & \frac{\partial F_1}{\partial x_2}(X^*)\\ 
	\frac{\partial F_1}{\partial x_2}(X^*) & \frac{\partial F_2}{\partial x_2}(X^*)
\end{pmatrix} \tag{3.1.5}\]
Then we can write (3.1.3) as follows
\[\label{e3.1.6}
\begin{cases}
dU^1_t = (a_1U^1_t + b_1U^2_t)dt + \sqrt{2T}dB^{H_1}_t\\
dU^2_t = (a_2U^1_t + b_2U^2_t)dt + \sqrt{2T}dB^{H_2}_t
\end{cases} \tag{3.1.6}
 \]
 
\subsection{Approximate model}
Now we can consider an approximate model corresponding to (3.1.4) or (3.1.6) by replacing $B^H_t = \begin{pmatrix} B^{H_1}_t \\B^{H_2}_t\end{pmatrix}$ by $B^{H,\epsilon}_t = \begin{pmatrix} B^{H_1,\epsilon_1}_t \\B^{H_2, \epsilon_2}_t\end{pmatrix}$ as introduced in section 1.
\[ \label{e3.2.1} dU^\epsilon_t = AU^\epsilon_tdt+ \sqrt{2T}dB^{H,\epsilon}_t \tag{3.2.1} \]
and
\[ \label{e3.2.2}\begin{cases}
	dU^{1,\epsilon_1}_t = (a_1U^{1,\epsilon_1}_t + b_1U^{2,\epsilon_2}_t)dt + \sqrt{2T}dB^{H_1,\epsilon_1}_t\\
	dU^{2,\epsilon_2}_t = (a_2U^{1,\epsilon_1}_t + b_2U^{2,\epsilon_2}_t)dt + \sqrt{2T}dB^{H_2,\epsilon_2}_t 
\end{cases} \tag{3.2.2}
\]
where $B^{H_j,\epsilon_j}_t \to B^{H_j}_t$ in the space $L^2  = L^2(\Omega)$ of square integrable processes, and
\[\label{e3.2.3} \norm{B^{H_j,\epsilon_j}_t - B^{H_j}_t}_{L^2} \le C_j(\alpha_j)\epsilon_j^{2H_j}, \alpha_j = H_j -1/2, j=1,2 \tag{3.2.3}\]
with $C(\alpha_j)$ depending only  on $\alpha_j$ (see [7]).
\subsection{Convergence}
Before finding the solution $U^\epsilon_t$ of (3.2.1) we will prove a convergence theorem.

 \begin{thm} The solution $U^\epsilon_t$ of (3.2.1) tend to the solution $U_t$ of (2.7) in $L^2$ as $\epsilon \to 0$, that means
 \[U^{1,\epsilon_1}_t \to U^1_t \text{ in } L^2 \text{ as } \epsilon_1 \to 0\]
 and
  \[U^{2,\epsilon_2}_t \to U^2_t \text{ in } L^2 \text{ as } \epsilon_2 \to 0\]
\end{thm}

\begin{pf} We have
\[ \label{e3.3.1}U^1_t  - U^{1,\epsilon_1}_t = \int_0^t \left[a_1(U^1_t-U^{1,\epsilon_1}) + b_1(U^2_t-U^{2,\epsilon_2}_t) \right]ds + \sqrt{2T}(B^{H_1}_t - B^{H_1,\epsilon_1}_t) \tag{3.3.1}\]
\[\label{e3.3.2} U^2_t  - U^{2,\epsilon_2}_t = \int_0^t \left[a_2(U^1_t-U^{1,\epsilon_1}) + b_2(U^2_t-U^{2,\epsilon_2}_t) \right]ds + \sqrt{2T}(B^{H_2}_t - B^{H_2,\epsilon_2}_t) \tag{3.3.2}\]
\end{pf}

Put $M = \max\left( \frac{\abs{a_1}}{2},\frac{\abs{a_2}}{2}, \frac{\abs{b_1}}{2},\frac{\abs{b_2}}{2} \right)$. Then it follows from (3.3.1), (3.3.2) and (2.5) that 
\[\label{e3.3.3}
\norm{U^1_t  - U^{1,\epsilon_1}_t } + \norm{U^2_t  - U^{2,\epsilon_2}_t } \le C(\alpha,\epsilon) + M\int_0^1 \left[ \norm{U^1_s-U^{1,\epsilon_1}_s} + \norm{U^2_t-U^{2,\epsilon_2}_t} \right]ds \tag{3.3.3}
\]
where
\[C(\alpha, \epsilon) = \sqrt{2T} \left(C_1(\alpha_1)\epsilon_1^{1+2\alpha_1} + C_2(\alpha_2)\epsilon_2^{1+2\alpha_2} \right) \to 0\]
as $\epsilon = (\epsilon_1, \epsilon_2) \to 0$
Then an application of Gronwall inequality gives us
\[\norm{U^1_t  - U^{1,\epsilon_1}_t } + \norm{U^2_t  - U^{2,\epsilon_2}_t }\le C(\alpha,\epsilon) \exp{(Mt)} \to 0 \]
as $\epsilon \to 0$. So that
\[\label{e3.3.4}U^\epsilon_t = \begin{pmatrix}U^{1,\epsilon_1}_t \\U^{2,\epsilon_2}_t\end{pmatrix} \to \begin{pmatrix}U^1_t \\U^2_t\end{pmatrix} = U_t \text{ in } L^2 \text{ as } \epsilon =(\epsilon_1,\epsilon_2) \to 0 \tag{3.3.4}\]
\subsection{Solution $U^\epsilon_t$ of (3.2.1)}

Now we have to  find the solution of $U^\epsilon_t$ in its explicit form. We have
\[dU^\epsilon_t = AU^\epsilon_tdt + \sqrt{2T}dB^{H,\epsilon}_t\]
where
\[A = \begin{pmatrix} a_1& b_1\\ a_2&b_2 \end{pmatrix}, U^\epsilon_t =  \begin{pmatrix} U^{H_1,\epsilon_1}_t\\ U^{H_2,\epsilon_2}_t\end{pmatrix}, B^\epsilon_t =  \begin{pmatrix} B^{H_1,\epsilon_1}_t\\ B^{H_2,\epsilon_2}_t\end{pmatrix}\]
It is easy to see that
\[\label{e3.4.1}U^\epsilon_t  = \int_0^te^{A(t-s)}\sqrt{2T}dB^{H,\epsilon}_t \tag{3.4.1}\]
A method by Oksendal can be applied ([6], page 99) to obtain
\[\label{e3.4.2}e^{A(t-s)} = \frac{e^{-\lambda(t-s)}}{\xi}[\xi\cos{\xi(t-s)} + \lambda \sin{\xi(t-s) }\begin{pmatrix}1&0\\0&1 \end{pmatrix} +A\sin{\xi(t-s)}] \tag{3.4.2}\]
where
\[\label{e3.4.3}\lambda = -\frac{b_2}{2}, \xi = \abs{a_2 - \frac{b^2_2}{4}}\tag{3.4.3} \]
(3.4.2) can be rewritten also as
\[ \label{e3.4.4} \begin{split}
e^{A(t-s)} =& \frac{e^{-\lambda(t-s)}}{\xi}\begin{bmatrix}\xi\cos{\xi(t-s)}+\lambda\sin{\xi(t-s)} \text{    }0\\ 0 \text{   } \xi\cos{\xi(t-s)}+\lambda\sin{\xi(t-s)}  \end{bmatrix} +
\begin{pmatrix}a_1&a_2\\b_1&b_2 \end{pmatrix}\sin{\xi(t-s)}\\
 = & \frac{e^{-\lambda(t-s)}}{\xi}\begin{pmatrix}A_1&B_1\\A_2&B_2 \end{pmatrix}
 \end{split} \tag{3.4.4}\]
 where
 \[\label{e3.4.5}\begin{split}
 A_1 & =\xi\cos{(t-s)}+(\lambda+a_1)\sin{(t-s)} \\
 B_1 & =b_1\sin{(t-s)}\\
 A_2 &=a_2\sin(t-s)\\
 B_2 & = \xi\cos(t-s)+(\lambda+b_2)\sin(t-s)
 \end{split} \tag{3.4.5}\]
 Then it follow from (3.4.1), (3.4.3) and (3.4.4) that we can state now:
 \begin{thm}
 The solution of the approximate model (3.2.1) given by
 \[ \label{e3.4.6}U^\epsilon_t = \begin{pmatrix}U^{1,\epsilon_1}_t \\ U^{2,\epsilon_2}_t\end{pmatrix} = \frac{\sqrt{2T}}{\xi}\int_0^t\begin{pmatrix}A_1dB^{H_1,\epsilon_1}_s+B_1dB^{H_2,\epsilon_2}\\ A_2dB^{H_1,\epsilon_1}_s+B_2dB^{H_2,\epsilon_2}\end{pmatrix} \tag{3.4.6}\]
 where $A_1, B_1,  A_2, B_2$ and $\xi$ are defined in (3.4.5), (3.4.3) and (3.1.5).
 \end{thm}
 
 In account of the Convergence Theorem 3.1 we can obtain
 \begin{thm}
Given a steady state $X^*$ in the full fractional annealing range the state $X_t$ of the system is the limit in the spase $L^2(\Omega)$ of $X^\epsilon_t$ when $\epsilon$ tends to 0: 
 \[X_t = L^2 - \lim_{\epsilon \to 0}X^\epsilon_t \]
 where $X^\epsilon_t = X^* + U^\epsilon_t$ and $U^\epsilon_t$ is obtained from (3.4.6).
 \end{thm}
 
 \section*{Conclusion}
 Then each state $X_t$ of the annealing system can be considered as a $L^2-$limit of $(X^* +U^{\epsilon}_t)$ as $\epsilon \to 0$ where $X^*$ is a steady state of full annealing range.
 \section*{Acknowledgments}
 This research is funded by Vietnam National Foundation for Science and Technology Development (NAFOSTED) under Grant 101.02-2011.12.

\end{document}